\documentclass[11pt,reqno]{book}
\usepackage{fullpage,enumerate}
\usepackage{amsfonts}
\usepackage[all]{xy}
\usepackage{amssymb}
\usepackage{comment}
\usepackage{times}
\usepackage{graphicx}
\usepackage{amsmath,bm}
\usepackage[hidelinks]{hyperref}
\usepackage{tikz-cd}
\usepackage{tikz}
\usepackage{appendix}
\usepackage{color}\usepackage{float}
\usepackage[margin=1.3in,marginparwidth=1.1in,marginparsep=0.1in]{geometry}
\usepackage{marginnote,color}
\usepackage{import}
\usepackage[normalem]{ulem}
\usepackage{mathrsfs}
\usepackage{epsfig,amscd,amsthm,mathtools,fourier}
\usepackage{wrapfig}

\newtheorem{theorem}{Theorem}[section]
\newtheorem{corollary}[theorem]{Corollary}
\newtheorem{lemma}[theorem]{Lemma}
\newtheorem{proposition}[theorem]{Proposition}

\theoremstyle{definition}
\newtheorem{definition}[theorem]{Definition}
\newtheorem{example}{Example}[theorem]

\newtheorem{remark}{Remark}[theorem]
\newtheorem{exercise}{Exercise}[subsection]

\newcommand{\Z}{{\mathbb Z}}
\newcommand{\R}{{\mathbb R}}
\newcommand{\C}{{\mathbb C}}

\renewcommand{\P}{{\mathbb P}}

\newcommand{\G}{{\mathbb G}}

\newcommand{\kl}{{\mathcal L}}
\newcommand{\km}{{\mathcal M}}

\newcommand{\ko}{{\mathcal O}}

\DeclareMathAlphabet{\mathsc}{U}{rsfs}{m}{n}

\newtheorem{lemma and definition}[theorem]{Lemma and Definition}

\newcommand{\bp}{\boldsymbol{p}}

\newcommand{\bq}{\boldsymbol{q}}

\DeclareMathOperator{\ord}{ord}

\DeclareMathOperator{\rank}{rank}

\DeclareMathOperator{\Spec}{Spec}

{\end{list}}

\newcounter{Enumi}
{\end{list}}

\begin{document}
\author{Ilya Tyomkin}
\title{The Geometry of Severi Varieties\footnote{This is an appendix to the second edition of the book {\em Introduction to Singularities and Deformations} by Greuel, Lossen, and Shustin \cite{GLS07}}}

\maketitle
\setcounter{chapter}{3}
\setcounter{section}{6}

In this appendix, we summarize known results on the geometry of Severi varieties on toric surfaces – the varieties parameterizing integral curves of a given geometric genus in a given linear system. Till the last decade, Severi varieties were studied exclusively in characteristic zero. In particular, in the ’80-s, Zariski proved that a general plane curve of a given genus is necessarily nodal and gave a dimension-theoretic characterization of the Severi varieties \cite{Zar82t}. A few years later, Harris showed that the classical Severi varieties are irreducible \cite{Har86t}. The geometry of Severi varieties is much subtler on general toric surfaces, especially in positive characteristic. In the appendix, we discuss in particular recent examples of reducible Severi varieties \cite{Tyo14, LT23} and of components of Severi varieties parameterizing non-nodal curves in positive characteristic \cite{Tyo13, Tyo14}. We explain the new tools coming from tropical geometry that allowed us to generalize the theorems of Zariski and Harris to arbitrary characteristic in the classical case of curves on the projective plane \cite{CHT23}. Finally, we discuss the results about the adjacency of Severi varieties for different genera following \cite{CHT22t}.

\subsection{Introduction}

In 1921, Francesco Severi introduced what nowadays is called the classical {\em Severi varieties} to be the loci of reduced plane curves of degree $d$ that have $\delta$ nodes as their only singularities, \cite{Severi}. It is easy to check that the classical Severi varieties $V^{d,\delta}\subseteq |\ko_{\P^2}(d)|$ are locally closed subsets, hence admitting a natural structure of quasi-projective varieties.

The original motivation of Severi for considering these varieties was to provide an algebraic proof of the irreducibility of the moduli spaces $\km_g$ of smooth projective curves (or compact Riemann surfaces) of genus $g$. The connectedness of those spaces was already known to Felix Klein in 1882, who deduced it from the connectedness property of the space of genus-$g$ degree-$d$ coverings of the projective line due to earlier works of Lüroth and Clebsch \cite{Lur71, Cle73}.

The idea of Severi was as follows. By the Riemann-Roch theorem, any irreducible smooth projective curve $C$ of genus $g$ admits an embedding into a projective space as a curve of large enough degree, e.g., $d=2g+1$. Projecting this curve onto a general plane gives rise to a nodal plane curve of degree $d$ with $\delta:=\binom{d-1}{2}-g$ nodes, whose normalization is $C$, see e.g., \cite[Corollary~IV.3.11]{Har77}. Therefore, there is a natural surjective map from the union of components $V^{d,\delta,{\rm irr}}\subseteq V^{d,\delta}$ parametrizing irreducible curves to the moduli space $\km_g$, which reduces the irreducibility problem of the moduli space $\km_g$ to the Severi problem -- the question whether the Severi varieties $V^{d,\delta,{\rm irr}}$ are irreducible. 

Severi showed that for any degree $d\ge 1$ and any $0\le \delta\le \binom{d}{2}$, the variety $V^{d,\delta}\subseteq |\ko_{\P^2}(d)|$ is equidimensional of codimension $\delta$. Thus, its dimension is given by $\dim(V^{d,\delta})=3d+g-1$ since $\dim(|\ko_{\P^2}(d)|)=\binom{d+2}{2}-1$ and $\delta=\binom{d-1}{2}-g$. Severi showed that the varieties $V^{d,\delta}$ are smooth and that $V^{d,\delta'}\subseteq \overline{V^{d,\delta}}$ for any $\delta'\ge \delta$. Furthermore, he described the branches of $V^{d,\delta}$ and $V^{d,\delta,{\rm irr}}$ along $V^{d,\delta'}$.

Severi's approach to the irreducibility of $V^{d,\delta,{\rm irr}}$ was the following. Set $\delta':=\binom{d}{2}$. First, show that for any irreducible component $V\subseteq V^{d,\delta}$, the closure $\overline{V}$ contains the Severi variety $V^{d,\delta'}$ parametrizing unions of $d$ lines, and second, use the description of the branches of $V^{d,\delta,{\rm irr}}$ along $V^{d,\delta'}$ to show that they all belong to the same component. However, Severi's proof of the first step contained a gap, and Severi problem remained open for more than 60 years; see \cite{Zar82t, Ful82} for more details.

Severi problem was settled in the case of characteristic zero by Harris in 1986, \cite{Har86t}, and in arbitrary characteristic by Christ, He, and Tyomkin only in 2023, \cite{CHT23}. The first major step towards the resolution of Severi's problem was made by Zariski in 1982, who proved that in characteristic zero, if $W\subset |\ko_{\P^2}(d)|$ is an irreducible closed subvariety such that the general curve $[C]\in W$ is reduced and has geometric genus $g$, then $\dim(W)\le 3d+g-1$. Furthermore, the equality holds if and only if $C$ is nodal, \cite{Zar82t}. A closely related result was obtained independently by Arbarello and Cornalba, \cite{AC80}. Zariski's theorem played an important role in Harris' solution of the Severi problem.

\bgroup
\def\leftmark{3 Singularities in Arbitrary Characteristics}

Following the result of Zariski, it makes sense to consider the slightly bigger loci of curves of a given geometric genus without specifying the types of their singularities. Such loci are also called Severi varieties and can be defined for arbitrary linear systems on algebraic surfaces.

\begin{definition}
Let $X$ be an irreducible projective surface, $\kl$ a line bundle, and $g$ an integer. By the {\em Severi variety of genus-$g$ curves in the linear system $|\kl|$}, we mean the locus $V_{g,\kl}$ of reduced curves $C$ of geometric genus $g$ contained in the smooth locus of $X$ and such that no component of $C$ has negative self-intersection. The locus of integral curves in $V_{g,\kl}$ is denoted by $V^{\rm irr}_{g,\kl}$, and is called the {\em Severi variety of genus-$g$ integral curves in the linear system $|\kl|$.}
\end{definition}

One can show that these more general Severi varieties are locally closed subsets of $|\kl|$, and therefore they are quasi-projective varieties, see, e.g., \cite[Lemma~2.6]{CHT23}. They have been intensively studied over the years for various surfaces, and a lot is known about their geometry, see e.g., \cite{Har86t, DH88b, CH98-2, CC99, Tyo07, Tes09, Tyo14, CFGK17, Z22, CGY23, CDGK23}. 

Till the last decade, Severi varieties were studied exclusively in characteristic zero. It turns out that in positive characteristic, their geometry and the geometry of a general curve of a given genus in a given linear system are subtler, and new phenomena occur. In particular, general curves are not necessarily nodal but may have other types of singularities depending on the characteristic, \cite{Tyo13}. 

In this appendix, we summarize known results about the geometry of Severi varieties, emphasizing the case of toric surfaces. We explain the new tools coming from tropical geometry that allow one to generalize the theorems of Zariski and Harris to arbitrary characteristic in the classical case of curves on the projective plane, \cite{CHT23}. 

We assume that the reader is familiar with the basics of toric geometry, in particular with constructions of toric varieties associated to lattice polytopes, the cone-to-orbit and face-to-orbit correspondences, etc. We refer the reader to one of the following introductory sources \cite{Dan78,Ful93,CLS11} for the background on toric geometry.

\subsubsection{Notation}

Throughout the appendix, $X$ denotes a projective surface, $\kl$ a line bundle, $g$ an integer, and $\delta$ the total $\delta$-invariant, i.e., the difference between the arithmetic genus of a general curve $C$ in the linear system $|\kl|$ and $g$, i.e.,
$\delta:=p_a(C)-g=1-\chi(\ko_X)+\chi(\kl^{-1})-g.$

By $M$ and $N$ we denote a pair of dual lattices, i.e., free abelian groups of finite rank. The lattice $M$ is the lattice of characters (or, equivalently, of monomials) on the algebraic torus $T:=\Spec K[M]$. We denote by $N_\R:=N\otimes_\Z\R$ and $M_\R:=M\otimes_\Z\R$ the real vector spaces spanned by $N$ and $M$, and by $\Delta\subset M_\R$ lattice polytopes, i.e., polytopes with vertices in $M$. 

For a lattice polytope $\Delta$, we denote by $\partial\Delta$ its boundary, by $\partial_i\Delta$ its faces of codimension-one, by $\Delta^\circ$ its interior, by $X_\Delta$ the associated toric variety, and by $\ko(\Delta)$ the tautological line bundle on $X_\Delta$. The monomial functions are denoted by $x^m$ for $m\in M$. Finally, we denote by $O_i\subset X_\Delta$ the codimension-one orbit corresponding to $\partial_i\Delta$ under the face-to-orbit correspondence, and by $\partial_iX_\Delta$ the closure of $O_i$ in $X_\Delta$. We denote the primitive inner normal to $\partial_i\Delta$ by $n_i$, i.e., $n_i$ is the primitive integral vector in the ray of the dual fan of $X_\Delta$ corresponding to the orbit $O_i$. An example to keep in mind is $M=N=\Z^2$ with the usual dot product, and $\Delta$ being the triangle with vertices $(0,0),(d,0),(0,d)$. Then $(X_\Delta,\ko(\Delta))=(\P^2, \ko_{\P^2}(d))$ and the sides of $\Delta$ correspond to the (one-dimensional torus orbits in the) three coordinate axis in $\P^2$.

We set $\partial X_\Delta:=\sum_i\partial_iX_\Delta$, and recall that for a toric variety $X_\Delta$, the canonical divisor $K_{X_\Delta}$ is equivalent to $-\partial X_\Delta$. If $\Delta$ is a polygon, then the dimension of the linear system $|\ko(\Delta)|$ is given by $\dim(|\ko(\Delta)|)=|\Delta\cap M|-1$, the arithmetic genus of a reduced curve $[C]\in|\ko(\Delta)|$ is given by the number of the integral inner points $|\Delta^\circ\cap M|$, and $K_{X_\Delta}\cdot C=-|\partial\Delta\cap M|$. Finally, if $(X,\kl)=(X_\Delta,\ko(\Delta))$ then we use the notation $V_{g,\Delta}$ and $V_{g,\Delta}^{\rm irr}$ for the Severi varieties $V_{g,\kl}$ and $V^{\rm irr}_{g,\kl}$, respectively.

\subsection{The Dimension}

We start with the classical case of the ground field $K$ of characteristic zero. In general, Severi varieties may have components of different dimensions; see \cite{CC99}. However, for rational surfaces, they are often equidimensional of codimension $\delta$. The standard approaches to studying the dimension of a Severi variety involve deformation theory, either of embedded curves \cite{Zar82t, KST13} or of parametrized curves \cite{AC81,AC81it,Har86t,Vak00,Tyo07}. The following theorem due to Kleiman and Shende is, perhaps, the most general result of this kind known to date; see \cite[Theorem~1]{KST13}.

\begin{theorem}\label{thm:dim}
    Assume $X$ is a rational surface over an algebraically closed field of characteristic zero. Let $V$ be a closed subset of $|\kl|$ that contains every $[C]\in |\kl|$ such that either
    \begin{enumerate}
        \item $C$ has a component $C_1$ with $-K_X\cdot C_1\le 0$, or
        \item $C$ has a non-immersed\footnote{By {\em non-immersed}, we mean that the differential of the normalization map $C_1^\nu\to C_1$ vanishes at some point.} component $C_1$ with $-K_X\cdot C_1= 1$.
    \end{enumerate} 
    Then $V_{g,\kl}\setminus V$ has codimension $\delta$ in $|\kl|$ at all its points (if any). Furthermore, its sublocus of immersed curves is open and dense, and is smooth off $V$.
\end{theorem}

\begin{remark}\label{rem:dim}
    Assume that $V_{g,\kl}\setminus V$ is not empty.  Then there is a reduced curve $[C]\in |\kl|$ all of whose components have non-negative self-intersection. Thus, by \cite[Lemma~7]{KST13}, $\dim|\kl|=\frac{1}{2}C\cdot(C-K_X)$, and therefore $V_{g,\kl}\setminus V$ is equidimensional of dimension $-K_X\cdot C+g-1$ by Theorem~\ref{thm:dim} and the adjunction formula. Indeed,
    \begin{equation*}
        \begin{split}
            \dim\left(V_{g,\kl}\setminus V\right)={} & \dim|\kl|-\delta =\\
            &\frac{1}{2}C\cdot(C-K_X)-p_a(C)+g =\\ 
            &\frac{1}{2}C\cdot(C-K_X)-\frac{1}{2}C\cdot(C+K_X)-1+g =\\
            &-K_X\cdot C+g-1.
        \end{split}
    \end{equation*}
\end{remark}

\begin{theorem}\label{thm:dimtor}
    Let $(X,\kl)=(X_\Delta,\ko(\Delta))$ be the polarized toric surface associated to a lattice polygon $\Delta$ over an algebraically closed field of characteristic zero, and $0\le g\le |\Delta^\circ\cap M|$ an integer. Then $V^{\rm irr}_{g,\Delta}\ne\emptyset$ and $V_{g,\Delta}$ is equidimensional of dimension $\dim(V_{g,\Delta})=|\partial\Delta\cap M|+g-1.$
\end{theorem}

\begin{proof}
    Let us first show that $V^{\rm irr}_{0,\Delta}\ne\emptyset$. As usual, let $n_i\in N$ be the primitive inner normal vector to the side $\partial_i\Delta$, and denote the length of $\partial_i\Delta$ by $k_i$. Set $k:=\sum k_i$ and pick $k$ points $\{p_{i,j}\}_{1\le j\le k_i}$ in general position in $\P^1$. Consider the homomorphism $M\to K(t)$ given by 
    $$m\mapsto \prod (t-p_{i,j})^{(n_i,m)}.$$
    It induces a rational map $\P^1\dashrightarrow T$ defined on the complement of $\{p_{i,j}\}$. We denote by $f$ the composition of this map with the inclusion $T\subset X$. If $O_i$ is the torus orbit corresponding to the side $\partial_i\Delta$, then the affine toric subvariety $T\cup O_i$ is the spectrum of the monoid algebra $K[M_i]$, where $M_i=\{m\in M\,|\, (n_i,m)\ge 0\}$. Therefore, $f$ extends to $p_{i,j}$, and $f(p_{i,j})\in O_i$. Moreover, $f^*(\partial_i X)=\sum_j p_{i,j}$, since ${\rm ord}_{\partial_i X}(x^m)=(n_i,m)$.

    We claim that $f\colon \P^1\to X_\Delta$ is birational onto its image. Since $f^*(\partial_i X)=\sum_j p_{i,j}$, it is sufficient to show that $f$ separates the points $\{p_{1,j}\}_j$. Pick a non-zero $m\in M$ such that $(n_1,m)=0$. Then the pullback $f^*(x^m)=\prod (t-p_{i,j})^{(n_i,m)}$ is independent of $\{p_{1,j}\}_j$, and since $\{p_{i,j}\}$ are general, $f^*(x^m)$ separates the points $\{p_{1,j}\}_j$. Hence so does $f$. Furthermore, we have shown that $f(\P^1)$ intersects each $\partial_iX$ transversally at $k_i$ points. We conclude that $f(\P^1)$ is an integral rational curve in the linear system $|\ko(\Delta)|$, which proves that $V^{\rm irr}_{0,\Delta}\ne\emptyset$.

    \begin{lemma}\label{lem:nodalityrat}
        For a general $[C]\in V^{\rm irr}_{0,\Delta}$, the rational curve $C$ is nodal.
    \end{lemma}
    
    Pick a general $[C]\in V^{\rm irr}_{0,\Delta}$. By Theorem~\ref{thm:dim}, the germ of the Severi variety at $C$ is smooth of codimension $\delta(C)$ in $|\ko(\Delta)|$. Furthermore, since $C$ is nodal by Lemma~\ref{lem:nodalityrat}, this germ is nothing but the space of equigeneric deformations of $C$. Therefore, the nodes of $C$ can be deformed independently, and thus $V^{\rm irr}_{g,\Delta}\ne\emptyset$ for any $0\le g\le |\Delta^\circ\cap M|$.

    Next, let us show that $V_{g,\Delta}\cap V=\emptyset$, where $V$ is as in Theorem~\ref{thm:dim}. Pick $[C]\in V_{g,\Delta}$, and let $C_1$ be an irreducible component of $C$. If $C_1\subset \partial X$, then it intersects exactly two other components of the boundary divisor $\partial X$. Furthermore, the intersection is transversal since $C$ is disjoint from the singular locus of $X$. Therefore, $-K_X\cdot C_1=2+C_1^2\ge 2$. Assume now that $C_1\not\subset \partial X$. Then, again, $-K_X\cdot C_1\ge 2$ since $C_1$ intersects $\partial X=-K_X$ at least twice. Indeed, any point of the toric surface admits an equivariant open affine neighborhood. Thus, if $C_1$ intersected $\partial X$ at a single point, it would be contained in an affine variety, which is absurd because $C_1$ is projective. We conclude that $[C]\notin V$, and hence $V_{g,\Delta}\cap V=\emptyset$.

    Now, by Theorem~\ref{thm:dim} and Remark~\ref{rem:dim}, the Severi variety $V_{g,\Delta}=V_{g,\Delta}\setminus V$ is equidimensional of dimension $-K_X\cdot C+g-1$. Finally, since $K_X\cdot C=-|\partial\Delta\cap M|$, the asserted dimension formula follows.
\end{proof}

{\renewcommand{\proofname}{Proof of Lemma~\ref{lem:nodalityrat}}
\begin{proof}
    Since $[C]\in |\ko(\Delta)|$, we have $-K_X\cdot C=|\partial\Delta\cap M|\ge 3$. If the inequality is strict, then the assertion follows from \cite[Proposition~2]{KST13}; otherwise the nodality can be verified by an explicit computation, see \cite[Lemma~3.5]{Shu05}.
\end{proof}
}

Unfortunately, the deformation-theoretic techniques used to prove Theorem~\ref{thm:dim} fail in the case of positive characteristic. Nevertheless, at least in the case of toric surfaces, one can use tropical tools to prove the following result, see \cite[Proposition~2.7]{CHT23} and \cite[Theorem~1.2]{Tyo13}.

\begin{theorem}\label{thm:dimest}
    Let $(X,\kl)=(X_\Delta,\ko(\Delta))$ be the polarized toric surface associated to a lattice polygon $\Delta$ over an algebraically closed field of arbitrary characteristic, and let $g$ be an integer. Then  $V_{g,\Delta}$ is either empty or equidimensional of dimension $\dim(V_{g,\Delta})=|\partial\Delta\cap M|+g-1.$
\end{theorem}

In section~\ref{sect:proofdimest}, we will explain how to use tropical techniques to prove the upper bound on the dimension: $\dim(V_{g,\Delta})\le|\partial\Delta\cap M|+g-1$.
\begin{remark}
    In the classical case $(X_\Delta,\ko(\Delta))=(\P^2,\ko_{\P^2}(d))$, one can prove more. Namely, by \cite[Theorem~6.1]{CHT23}, for any $1-d\le g\le \binom{d-1}{2}$, the Severi variety $V_{g,d}$ is non-empty of pure dimension $3d+g-1$ in any characteristic. Furthermore, if $0\le g\le \binom{d-1}{2}$ then $V_{g,d}^{\rm irr}\ne\emptyset$.
\end{remark}

\subsection{The Geometry of a General Curve of a given Genus}
In this subsection, we address the question about the geometry of a general curve $[C]\in V_{g,\kl}$. This question was first investigated by Zariski in \cite[Theorem~2]{Zar82t}, who proved the following:
\begin{theorem}[Zariski's Theorem]\label{thm:Zarthm}
    Assume that the ground field $K$ has characteristic zero, and let $W\subset |\ko_{\P^2}(d)|$ be an irreducible closed subvariety. Assume that the general curve $[C]\in W$ is reduced and has geometric genus $g$. Then $\dim(W)\le 3d+g-1$. Furthermore, the equality holds if and only if $C$ is nodal and $W$ is the closure of a component of $V_{g,d}$.
\end{theorem}

Several generalizations of this result can be found in the literature; see, e.g., \cite[Theorem~2.8]{Tyo07} and \cite[Proposition~2]{KST13}. In the case of Severi varieties on toric surfaces, the following is true.

\begin{theorem}\label{thm:gencurve}
    Let $(X,\kl)=(X_\Delta,\ko(\Delta))$ be the polarized toric surface associated to a lattice polygon $\Delta$ over an algebraically closed field of characteristic zero. Assume that $|\partial\Delta\cap M|>3$, and let $D\subset X$ be any curve. Then a general curve $[C]\in V^{\rm irr}_{g,\Delta}$ is nodal. Furthermore, $C$ intersects $D$ transversally; in particular, the intersection points are smooth points of $C$ and of $D$.
\end{theorem}
\begin{remark}
    Theorem~\ref{thm:gencurve} implies that the system $V^{\rm irr}_{g,\Delta}$ has no base points, and in particular, for a general $[C]\in V^{\rm irr}_{g,\Delta}$, the curve $C$ does not pass through the zero-dimensional orbits of $X$. Indeed, to see this, it suffices to pick a curve $D$ having a singularity at a given point.
\end{remark}
\begin{proof}
    Since $[C]\in \ko(\Delta)$, we have $-K_X\cdot C=|\partial\Delta\cap M|>3$. The nodality assertion now follows from \cite[Proposition~2]{KST13} and the transversality from \cite[Proposition~17~(4)]{KST13}.
\end{proof}

\begin{remark}
    We expect the theorem to hold true even without the assumption $|\partial\Delta\cap M|>3$. If $|\partial\Delta\cap M|\le 3$, then the equality holds and $-K_X\cdot C=3$. By \cite[Proposition~2]{KST13}, if $[C]\in V^{\rm irr}_{g,\Delta}$ is general, then it is necessarily immersed (in the topological sense), and no two branches are tangent to each other. It remains to rule out the possibility for simple multiple points other than nodes. If $g=0$, this was verified by Shustin \cite[Lemma~3.5]{Shu05}, but the case of positive genus requires an additional argument. By \cite[Proposition~17~(4)]{KST13}, $C$ intersects any given curve transversally. Finally, since the irreducible components of curves parametrized by the Severi variety $V_{g,\Delta}$ vary independently, and any singular curve $C_1$ on a toric surface satisfies $-K_X\cdot C_1\ge 3$, the only possible obstacle for the nodality of a general curve $[C]\in V_{g,\Delta}$ is the existence of simple multiple points of its irreducible components.
\end{remark}

Next, let us discuss the case of positive characteristic. It turns out that Theorem~\ref{thm:gencurve} completely fails in this case. Neither a general curve $[C]\in V^{\rm irr}_{g,\Delta}$ must be nodal, nor it must intersect a given curve transversally. Below, we give the simplest example of such behavior. More examples and details can be found in \cite[\S4.1]{Tyo13}.

Consider the lattice triangle $\Delta$ with vertices $(0,0), (2,1), (1,2)$, and the corresponding Severi variety $V^{\rm irr}_{0,\Delta}$. Notice that the three zero-dimensional orbits are singular points of $X$, and therefore, any $[C]\in V^{\rm irr}_{0,\Delta}$ contains none of them. As usual, denote the primitive inner normals to $\partial\Delta$ by $n_1:=(-1,2), n_2:=(-1,-1), n_3:=(2,-1)$, and let $O_i\subset X$ be the corresponding one-dimensional orbits. Let $[C]\in V^{\rm irr}_{0,\Delta}$ be any curve, and $f\colon \P^1\to C\subset X_\Delta$ its normalization. Since the integral lengths of the sides of $\Delta$ are $1$, exactly three points are mapped to the boundary divisor. Without loss of generality, we may assume that $0$ is mapped to $O_1$, $1$ to $O_2$ and $\infty$ to $O_3$. Thus, $${\rm div}(f^*(x^m))=(n_1,m)\cdot 0+(n_2,m)\cdot 1+(n_3,m)\cdot\infty$$ for any $m\in \Z^2$, and therefore $f^*(x^m)=\chi(m)t^{(n_1,m)}(t-1)^{(n_2,m)}$ for some multiplicative character $\chi\colon \Z^2\to K^\times$. We claim that if the characteristic is $3$, then there is $t\in \P^1\setminus\{0,1,\infty\}$ such that the differential of $f$ vanishes at $t$. Indeed, the vanishing of (log-)derivatives of $f^*(x^m)$ at $t$ is equivalent to $$0=\frac{(n_1,m)}{t}+\frac{(n_2,m)}{t-1}=\left(\frac{(-1,2)}{t}+\frac{(2,-1)}{t-1}, m\right),$$
and for $t=-1$ we have $\frac{(-1,2)}{t}+\frac{(2,-1)}{t-1}=(0,0)$ since the characteristic is $3$. Finally, since the arithmetic genus is $|\Delta^\circ\cap\Z^2|=1$, the curve $C$ is necessarily cuspidal. Furthermore, a similar computation shows that if the characteristic is $3$, then any pair of curves $[C_1],[C_2]\in V^{\rm irr}_{0,\Delta}$ intersect at a single point, and the intersection index at this point is $3$.


Let us present another point of view on the example above that admits generalizations to higher genera and other toric surfaces. Consider the sublattice $M\subset\Z^2$ generated by $\partial\Delta\cap\Z^2$. Then $\Delta$ is integral with respect to $M$, and the corresponding polarized toric  surface $X_M$ is isomorphic to $(\P^2,\ko_{\P^2}(1))$. There is a natural toric morphism $\pi\colon X_M\to X$, and $X$ is the quotient of $X_M$ by a subgroup of the torus isomorphic to $\mu_3$. Furthermore, the lines in $X_M$ project to rational curves in the linear system $|\kl|$, and since $\dim(V^{\rm irr}_{0,\Delta})=2=\dim |\ko_{\P^2}(1)|$, the pushforward of the linear system of lines on $X_M$ to $X$ gives rise to an irreducible component of $V^{\rm irr}_{0,\Delta}$. In fact, it is the whole of $V^{\rm irr}_{0,\Delta}$, since the latter is irreducible. Now, if the characteristic is $3$, then the map $\pi$ is {\em bijective}. Hence all singularities of curves parametrized by $V^{\rm irr}_{0,\Delta}$ have to be unibranch. Furthermore, any pair of such curves must intersect at a single point.

\subsection{The (Ir)Reducibility}
\subsubsection{The Severi problem}
As explained in the introduction, the original motivation of Severi for considering Severi varieties was to provide an algebraic proof of the irreducibility of the moduli spaces $\km_g$ of smooth projective curves (or compact Riemann surfaces) of genus $g$. 
By the Riemann-Roch theorem, any smooth projective curve $C$ of genus $g$ admits an embedding into a projective space as a curve of degree $d=2g+1$. Projecting this curve onto a general plane gives rise to a nodal plane curve of degree $d$ whose normalization is $C$, see e.g., \cite[Corollary~IV.3.11]{Har77}. Therefore, there is a natural surjective rational map from the Severi variety $V^{\rm irr}_{g,d}$ to $\km_g$, which reduces the irreducibility problem of the moduli space $\km_g$ to the Severi problem -- the question whether the Severi varieties $V^{\rm irr}_{g,d}$ are irreducible. 
\begin{remark}
    To define a morphism from an open subset of $V^{\rm irr}_{g,d}$ to $\km_g$ one must equinormalize the tautological family of curves over it. Since the locus of nodal curves $V^{d,\delta,{\rm irr}}\subseteq V^{\rm irr}_{g,d}$ is smooth by Theorem~\ref{thm:dim}, see also \cite[Propoisition~6.3]{CHT23} for the positive characteristic case, the tautological family over $V^{d,\delta,{\rm irr}}$ is equinormalizable by \cite{CL06}, see also \cite[Theorem~33]{GP22}. Thus, the induced map $V^{d,\delta,{\rm irr}}\to \km_g$ is surjective.
\end{remark}
Severi provided an argument for the irreducibility of $V^{\rm irr}_{g,d}$ based on the degeneration of plane curves to (non-generic) unions of lines. However, his proof contained a gap, and Severi problem remained open for more than 60 years. It was settled in the case of characteristic zero by Harris in 1986, see \cite{Har86t}, and in arbitrary characteristic by Christ, He, and Tyomkin only in 2023, see \cite{CHT23}:

\begin{theorem}\label{thm:irr}
    Let $d\ge 1$ and $0\le g\le \frac{(d-1)(d-2)}{2}$ be integers. Then, the Severi variety $V^{\rm irr}_{g,d}$ is irreducible. 
\end{theorem}

\begin{remark}
    Using Theorem~\ref{thm:irr}, one obtains the first known proof of the irreducibility of the moduli spaces $\km_g$ in positive characteristic that involves no reduction to the characteristic zero case.  Deligne-Mumford's classical proof, \cite{DM69}, proceeds by first reducing to the characteristic zero case and then using Teichm\"uller theory over $\mathbb C$. As a result, it is based on transcendental methods. The first algebraic proof is due to Fulton, who replaced the transcendental methods of Teichm\"uller theory by the study of the Harris-Mumford's compactification of Hurwitz schemes \cite{HM82}.  
\end{remark}

\begin{corollary}\label{cor:nod}
    Let $[C]\in V^{\rm irr}_{g,d}$ be a general curve. Then $C$ is nodal.
\end{corollary}
\begin{proof}
    Since $V^{\rm irr}_{g,d}$ is irreducible and nodality is an open condition, it is sufficient to construct a nodal curve of degree $d$ and geometric genus $g$. To do so, consider a union $C_0$ of $d$ general lines. Let $L$ be one of those lines, and pick a collection $Q$ of $d+g-1$ nodes that contains the $d-1$ nodes of $C_0$ belonging to $L$. Set $P:=C_0^{\rm sing}\setminus Q$. It is enough to show that there exists a deformation of $C_0$ preserving the nodes in $P$ and smoothing out those in $Q$. Indeed, the general curve in such a deformation is a nodal curve of geometric genus $g$, and it is irreducible since $C_0\setminus P$ is connected. The existence of such a deformation is clear since the equigeneric locus of $C_0$ in $|\ko_{\P^2}(d)|$ consists of all unions of $d$ distinct lines and hence has codimension $$\frac{d(d+3)}{2}-2d=\frac{d(d-1)}{2}.$$ Thus, the $\frac{d(d-1)}{2}$ nodes of $C_0$ can be deformed independently.
\end{proof}

\begin{remark}
Corollary~\ref{cor:nod} is nothing but a part of Zariski's Theorem. In characteristic zero, it was proved in order to show the irreducibility of the Severi varieties, and Harris' solution of the Severi problem uses it heavily. Curiously enough, in positive characteristic, we do not know how to prove this corollary without proving the irreducibility of the Severi varieties first.
\end{remark}

\subsubsection{Examples of reducible Severi varieties}
Although for non-rational surfaces, Severi varieties may be reducible and even non-equidimensional, see \cite{CC99}, the first non-trivial examples of reducible Severi varieties on toric surfaces were discovered only recently, see \cite{Tyo13, Tyo14}. Let us describe these examples following \cite{LT23}.


\begin{theorem}\label{thm:general}
Let $\Delta\subset \R^2$ be a lattice polygon and $g\ge 1$ and integer. Then the number of irreducible components of the Severi variety $V_{g,\Delta}^{\rm irr}$ is bounded from below by the number of affine sublattices $M\subseteq\Z^2$ for which the following two conditions hold: 
\begin{enumerate}
    \item $\partial\Delta\cap M=\partial\Delta\cap\Z^2$, and
    \item $|\Delta^\circ\cap M|\ge g$.
\end{enumerate}
\end{theorem}

\begin{example}
    Let $\Delta\in\R^2$ be the triangle with vertices $(0,0), (1,0), (-16,105)$, and set $g:=1$. Then the sublatices $M$ satisfying the two conditions of the theorem are $\Z\oplus d\Z$ for all positive divisors of $105$ excluding $105$ itself. Thus, the number of the irreducible components of the Severi variety $V_{1,\Delta}^{\rm irr}$ is at least $\sigma(3\cdot5\cdot7)-1=2^3-1=7$. 
\end{example}

\setcounter{remark}{1}
\begin{remark}
    In genus zero, the Severi varieties $V_{0,\Delta}^{\rm irr}$ are always irreducible by \cite[Proposition~4.1]{Tyo07}. In genus one, a recent result of Barash and Tyomkin shows that the inequality in Theorem~\ref{thm:general} is, in fact, an equality, \cite{BT24}. Thus, seven is the actual number of the irreducible components in the example above. Finally, if the genus is greater than one, then the bound in Theorem~\ref{thm:general} is not always sharp, as the following theorem shows.
\end{remark}

\begin{theorem}\label{thm:kite}
Let $k',k$ be non-negative integers such that $k'\ge k$ and $k'>0$. Let $\Delta_{k,k'}$ be the polygon with vertices $(0,0),(\pm 1,k), (0,k+k')$, and $(X_\Delta,\ko(\Delta))$ the corresponding polarized toric surface. Then for any integer $g\geqslant 1$, the number of irreducible components of $V_{g,\Delta}^{\rm irr}$ is bounded from below by the number of affine sublattices $M\subseteq\Z^2$ satisfying the two conditions of Theorem~\ref{thm:general} but counted with the following multiplicities: if the index of $M$ in $\Z^2$ is odd, then the multiplicity is equal to the number of integers $0\le \kappa\le \min\{|\Delta^\circ\cap M|-g, g\}$ of the same parity as $|\Delta^\circ\cap M|-g$, and is one otherwise.
\end{theorem}

\setcounter{figure}{2}
\begin{figure}[ht]
\begin{center}
\tikzset{every picture/.style={line width=0.75pt}} 
\begin{tikzpicture}[x=0.75pt,y=0.75pt,yscale=-.35,xscale=.35]
\draw [color={rgb, 255:red, 0; green, 0; blue, 0 }  ,draw opacity=1 ][line width=0.7]  [dash pattern={on 1.69pt off 2.76pt}]  (80,270) .. controls (120.5,271) and (103.5,188) .. (130,270) ;
\draw [fill={rgb, 255:red, 228; green, 228; blue, 232 }  ,fill opacity=1 ]   (30,170.01) -- (80,20) -- (130,170.01) -- (80,270) -- cycle ;
\draw [shift={(30,170.01)}, rotate = 243.43] [color={rgb, 255:red, 0; green, 0; blue, 0 }  ][fill={rgb, 255:red, 0; green, 0; blue, 0 }  ][line width=0.75]      (0, 0) circle [x radius= 3.35, y radius= 3.35]   ;
\draw [shift={(30,170.01)}, rotate = 288.43] [color={rgb, 255:red, 0; green, 0; blue, 0 }  ][fill={rgb, 255:red, 0; green, 0; blue, 0 }  ][line width=0.75]      (0, 0) circle [x radius= 3.35, y radius= 3.35]   ;
\draw    (80,220) ;
\draw [shift={(80,220)}, rotate = 0] [color={rgb, 255:red, 0; green, 0; blue, 0 }  ][fill={rgb, 255:red, 0; green, 0; blue, 0 }  ][line width=0.75]      (0, 0) circle [x radius= 3.35, y radius= 3.35]   ;
\draw [shift={(80,220)}, rotate = 0] [color={rgb, 255:red, 0; green, 0; blue, 0 }  ][fill={rgb, 255:red, 0; green, 0; blue, 0 }  ][line width=0.75]      (0, 0) circle [x radius= 3.35, y radius= 3.35]   ;
\draw    (80,170) ;
\draw [shift={(80,170)}, rotate = 0] [color={rgb, 255:red, 0; green, 0; blue, 0 }  ][fill={rgb, 255:red, 0; green, 0; blue, 0 }  ][line width=0.75]      (0, 0) circle [x radius= 3.35, y radius= 3.35]   ;
\draw [shift={(80,170)}, rotate = 0] [color={rgb, 255:red, 0; green, 0; blue, 0 }  ][fill={rgb, 255:red, 0; green, 0; blue, 0 }  ][line width=0.75]      (0, 0) circle [x radius= 3.35, y radius= 3.35]   ;
\draw    (80,120) ;
\draw [shift={(80,120)}, rotate = 0] [color={rgb, 255:red, 0; green, 0; blue, 0 }  ][fill={rgb, 255:red, 0; green, 0; blue, 0 }  ][line width=0.75]      (0, 0) circle [x radius= 3.35, y radius= 3.35]   ;
\draw [shift={(80,120)}, rotate = 0] [color={rgb, 255:red, 0; green, 0; blue, 0 }  ][fill={rgb, 255:red, 0; green, 0; blue, 0 }  ][line width=0.75]      (0, 0) circle [x radius= 3.35, y radius= 3.35]   ;
\draw    (80,70) ;
\draw [shift={(80,70)}, rotate = 0] [color={rgb, 255:red, 0; green, 0; blue, 0 }  ][fill={rgb, 255:red, 0; green, 0; blue, 0 }  ][line width=0.75]      (0, 0) circle [x radius= 3.35, y radius= 3.35]   ;
\draw [shift={(80,70)}, rotate = 0] [color={rgb, 255:red, 0; green, 0; blue, 0 }  ][fill={rgb, 255:red, 0; green, 0; blue, 0 }  ][line width=0.75]      (0, 0) circle [x radius= 3.35, y radius= 3.35]   ;
\draw    (80,20) ;
\draw [shift={(80,20)}, rotate = 0] [color={rgb, 255:red, 0; green, 0; blue, 0 }  ][fill={rgb, 255:red, 0; green, 0; blue, 0 }  ][line width=0.75]      (0, 0) circle [x radius= 3.35, y radius= 3.35]   ;
\draw [shift={(80,20)}, rotate = 0] [color={rgb, 255:red, 0; green, 0; blue, 0 }  ][fill={rgb, 255:red, 0; green, 0; blue, 0 }  ][line width=0.75]      (0, 0) circle [x radius= 3.35, y radius= 3.35]   ;
\draw    (130,170.01) ;
\draw [shift={(130,170.01)}, rotate = 0] [color={rgb, 255:red, 0; green, 0; blue, 0 }  ][fill={rgb, 255:red, 0; green, 0; blue, 0 }  ][line width=0.75]      (0, 0) circle [x radius= 3.35, y radius= 3.35]   ;
\draw [shift={(130,170.01)}, rotate = 0] [color={rgb, 255:red, 0; green, 0; blue, 0 }  ][fill={rgb, 255:red, 0; green, 0; blue, 0 }  ][line width=0.75]      (0, 0) circle [x radius= 3.35, y radius= 3.35]   ;
\draw    (80,270) ;
\draw [shift={(80,270)}, rotate = 0] [color={rgb, 255:red, 0; green, 0; blue, 0 }  ][fill={rgb, 255:red, 0; green, 0; blue, 0 }  ][line width=0.75]      (0, 0) circle [x radius= 3.35, y radius= 3.35]   ;
\draw [shift={(80,270)}, rotate = 0] [color={rgb, 255:red, 0; green, 0; blue, 0 }  ][fill={rgb, 255:red, 0; green, 0; blue, 0 }  ][line width=0.75]      (0, 0) circle [x radius= 3.35, y radius= 3.35]   ;

\end{tikzpicture}
\caption{The kite $\Delta_{2,3}$.}
\label{fig:kite}
\end{center}
\end{figure}
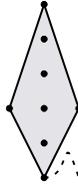

\begin{remark}
    The polygons $\Delta_{k,k'}$ are called {\em kites}. We shall mention, that even if the only sublattice satisfying conditions (a) and (b) is $M=\mathbb Z^2$, its multiplicity may be greater than one. For example, this is the case for the kite $\Delta_{2,3}$, and genus $g=2$.
\end{remark}

\subsection{The Degree}

Another fundamental problem regarding the geometry of the Severi varieties is the computation of their degrees. Several approaches to this problem have been developed by Kontsevich, Ran, Caporaso-Harris, Mikhalkin, and others. 

Let us mention that the degree of the Severi variety provides an answer to the following enumerative problem. Set $r:=\dim V_{g,\kl}^{\rm irr}$. {\em How many integral curves of the given degree and geometric genus that pass through $r$ points in general position?} Indeed, the locus of curves in the linear system $|\kl|$ passing through a given point is a hyperplane in the projective space $|\kl|$, and therefore, the degree of the Severi variety is equal to the number of points in its intersection with $r$ such hyperplanes, assuming that the intersection is transversal. The latter is not difficult to verify if a general curve $[C]\in V_{g,\kl}^{\rm irr}$ is nodal and the points are in a general position.

Let us start with the classical case of rational plane curves of degree $d$ in characteristic $0$. Set $N(d):=\deg(V^{\rm irr}_{0,d})$. Then there is a recursive formula due to Kontsevich \cite[Claim~5.2.1]{KM94}, which can be generalized to some Hirzebruch surfaces, see \cite{CH98-1}. 
\begin{theorem}[Kontsevich's formula]
    For $d\ge 2$, the following holds:
    $$N(d)=\sum_{k+l=d}N(k)N(l)k^2l\left[l\binom{3d-4}{3k-2}-k\binom{3d-4}{3k-1}\right].$$
\end{theorem}

To obtain recursive formulae for other surfaces or for curves of positive genus, one often has to consider the so-called {\em generalized Severi varieties} that parametrize irreducible curves of a given genus in a given linear system having prescribed tangencies to certain curves. In 1998, Caporaso and Harris proved a beautiful recursive formula for the degrees of the (generalized) Severi varieties in the planar case \cite[Theorem~1.1]{CH98-2}. To state their result we will need some preparations. 

Let $L\subset\P^2$ be a line, $\alpha,\beta$ two finite sequences of non-negative integers, and $\{p_{i,j}\}_{1\le j\le\alpha_i}\subset L$ a collection of points in general position. The {\em generalized Severi variety} $V^{d,\delta}(\alpha,\beta)$ is the locus of reduced nodal curves $C$ of degree $d$ with exactly $\delta$ nodes such that 
\begin{itemize}
    \item the line $L$ is not a component of $C$, and
    \item there exist points $q_{i,j},r_{i,j}$ in the normalization $C^\nu$ of the curve $C$ satisfying the following: $\nu^*(L)=\sum i(q_{i,j}+r_{i,j})$ and $\nu(q_{i,j})=p_{i,j}$ for all $i$ and $j$.
\end{itemize}

For two sequences $\epsilon$ and $\epsilon'$, we say that $\epsilon'\ge \epsilon$ if and only if $\epsilon'_i\ge \epsilon_i$ for all $i$, and in this case we set $\binom{\epsilon'}{\epsilon}:=\prod_i\binom{\epsilon'_i}{\epsilon_i}$. We denote by $e_k$ the zero sequence except for the value $1$ in the $k^{\rm th}$ place. Finally, set $I^\epsilon:=\prod_i i^{\epsilon_i}$

\begin{theorem}[Caporaso-Harris' formula]\label{thm:CH}
Assume that the characteristic of the ground field is zero, and set $N^{d,\delta}(\alpha,\beta):=\deg\left(V^{d,\delta}(\alpha,\beta)\right)$. Then,
    \begin{equation*}
        N^{d,\delta}(\alpha,\beta)=\sum_{k:\beta_k>0}k N^{d,\delta}(\alpha+e_k,\beta-e_k)+\sum I^{\beta'-\beta}\binom{\alpha}{\alpha'}\binom{\beta'}{\beta}N^{d-1,\delta'}(\alpha',\beta'),
    \end{equation*}
where the second sum is taken over all $\alpha',\beta',$ and $\delta'$ satisfying the following $\alpha'\le\alpha$, $\beta'\ge\beta$, $\delta'\le\delta$, and $\delta-\delta'+|\beta'-\beta|=d-1$. 
\end{theorem}

In the early 2000s, Mikhalkin developed a completely different novel approach to computing the degrees of (the generalized) Severi varieties. Unlike the formulae of Kontsevich and Caporaso-Harris, Mikhalkin's formula is not recursive and provides the answer in terms of enumeration of combinatorial objects, called {\em tropical curves}, with appropriate multiplicities. In the following section, we will define tropical curves, their degree, genus, and Mikhalkin's multiplicities. We will also describe the tropicalization procedure relating those to the algebraic curves and state Mikhalkin's Correspondence Theorem. Now, let us only state Mikhalkin's formula \cite[Theorem~1]{Mik05}.

\begin{theorem}[Mikhalkin's formula]\label{thm:MikFor}
  Let $(X_\Delta,\ko(\Delta))$ be the polarized toric surface associated to a lattice polygon $\Delta$ over an algebraically closed field of characteristic zero, $0\le g\le |\Delta^\circ\cap M|$ an integer, and $V^{\rm irr}_{g,\Delta}$ and $V_{g,\Delta}$ the corresponding Severi varieties. Let $\nabla$ be the reduced tropical degree dual to $\Delta$. Then the degree of $V_{g,\Delta}$ (resp., $V^{\rm irr}_{g,\Delta}$) is equal to the number of isomorphism classes of (resp., irreducible) plane tropical curves of degree $\nabla$ and genus $g$ passing through a tropically general collection of $r$ points and counted with Mikhalkin's multiplicities. 
\end{theorem}

\begin{remark}
    In 2007, Gathmann and Markwig considered the enumeration problem for plane tropical curves with prescribed tangencies to the line at infinity and proved that the resulting tropical invariants satisfy the Caporaso-Harris recursion \cite[Theorem~4.2]{GM07a}. Their proof is purely combinatorial. A version of Mikhalkin's Correspondence Theorem allows one to prove the equivalence of the tropical result and Theorem~\ref{thm:CH}. Therefore, one can deduce the algebraic Caporaso-Harris's recursion from the tropical one and vice versa. 
\end{remark}

\subsection{A brief Introduction to Tropical Curves and Tropicalization}

In this section, graphs are always finite graphs. They are allowed to have multiple edges, loops, and half-edges, which we call {\em legs}. If $\mathbb G$ is a graph, we denote by $V(\G),E(\G),$ and $L(\G)$ the sets of vertices, edges, and legs of $\G$, respectively. Set $\overline E(\mathbb G) := E(\G) \cup L(\G)$. For $e\in \overline E(\mathbb G)$, we use the notation $\vec e$ to indicate a choice of orientation on $e$. If $e\in L(\G)$ is a leg, then it will {\em always} be oriented away from the vertex. The edges will be considered with both possible orientations. We denote by ${\mathfrak h}(\vec e)$ and ${\mathfrak t}(\vec e)$ the head and the tail of an oriented edge $\vec e$. Similarly, ${\mathfrak t}(\vec l)$ denotes the tail of a leg $l$. For $v\in V(\mathbb G)$, we denote by ${\rm Star}(v)$ the {\em star} of $v$, i.e., the collection of oriented edges and legs having $v$ as their tail. In particular, ${\rm Star}(v)$ contains two oriented edges for every loop adjacent to $v$. The number of edges and legs in ${\rm Star}(v)$ is called the {\em valency} of $v$ and is denoted by ${\rm val}(v)$.

\subsubsection{Abstract and parametrized tropical curves}

\begin{definition}
    A {\em tropical curve} is a weighted metric graph $\Gamma=(\mathbb G,\ell)$ with ordered legs, i.e., $\mathbb G$ is a finite graph with ordered legs equipped with a {\em weight (or genus) function} $g\colon  V(\mathbb G)\rightarrow \mathbb Z_{\geq 0}$, and a {\em length function} $\ell\colon  E(\mathbb G)\rightarrow \mathbb R_{>0}$.
\end{definition}

The length function is usually extended to $\overline E(\mathbb G)$ by setting the length of any leg to be infinity. A tropical curve $\Gamma$ is called {\em irreducible} if the graph $\G$ is connected, and it is called {\em stable} if for every vertex $v\in V(\G)$ the following holds: $2g(v)-2+{\rm val}(v)\ge 1$. We define the {\em genus} of $\Gamma$ to be $g(\Gamma):=1-\chi(\mathbb G)+\sum_{v\in V(\mathbb G)}g(v)$, where as usual, $\chi(\G):=b_0(\G)-b_1(\G)$ denotes the Euler characteristic of $\G$.

\begin{remark}
    We will view tropical curves as polyhedral complexes by identifying the edges of $\G$ with bounded closed intervals of corresponding lengths in $\R$ and identifying the legs of $\G$ with semi-bounded closed intervals in $\R$.
\end{remark}


\begin{definition}
    A {\em parametrized tropical curve} is a pair $(\Gamma,h)$, where $\Gamma=(\mathbb G,\ell)$ is a tropical curve, and $h\colon  \Gamma\rightarrow N_\mathbb R$ is a map such that:
	
(a) for any $e\in \overline E(\mathbb G)$, the restriction $h|_e$ is an integral affine function, i.e., an affine function with an integral slope; and	
	
(b) for any vertex $v\in V(\mathbb G)$, the following {\em balancing condition} holds $$\sum_{\vec e\in{\rm Star}(v)}\frac{\partial h}{\partial \vec e}=0.$$
\end{definition}

\begin{remark}
    Note that the slope $\frac{\partial h}{\partial \vec e}\in N$ is not necessarily primitive, and its integral length, often called the {\it multiplicity} of $h$ along $e$, is the stretching factor of the affine map $h|_e$. In particular, the multiplicity vanishes if and only if $h$ contracts $e$.
\end{remark}

\begin{exercise}
    Let $h\colon \Gamma\to N_\R$ be a parametrized tropical curve. Show that $$\sum_{l\in L(\G)}\frac{\partial h}{\partial \vec l}=0.$$
\end{exercise}

We define the {\em extended degree} $\overline\nabla$ to be the sequence of slopes of $h$ along the legs of $\G$, and the {\em degree} $\nabla$ to be the subsequence of the non-zero slopes. We say that the degree $\nabla$ is {\em reduced} if it consists of primitive integral vectors. If the lattice $N$ has rank two and $\Delta$ is a lattice polygon in $M_\R$, then we say that the degree $\nabla$ is {\em dual} to $\Delta$ if $\nabla$ consists of outer normal vectors to the sides of $\Delta$ and for each side $\partial_j\Delta$, the sum of the multiplicities of the outer normals to $\partial_j\Delta$ in $\nabla$ is precisely the integral length $|\partial_j\Delta|$.

A parametrized tropical curve $(\Gamma, h)$ is called {\em stable} if so is $\Gamma$. An {\em isomorphism} of parametrized tropical curves $(\Gamma, h)$ and $(\Gamma', h')$ is an isomorphism of weighted metric graphs $\varphi\colon  \Gamma \to \Gamma'$ such that $h = h' \circ \varphi$. Finally, a {\em plane tropical curve} is a parametrized tropical curve for which $\rank(N)=2$.

\subsubsection{Moduli of parametrized tropical curves}
In this subsection we define the moduli spaces $M_{g, r, \nabla}^{\rm trop}$ of parametrized tropical curves of a given degree $\nabla$ with exactly $r$ contracted legs and genus $g$.  We assume that the contracted legs are $l_1,\dotsc,l_r$. Such parameter spaces can be considered either as polyhedral complexes (with an integral affine structure) as in \cite{GM07}, or, following Abramovich, Caporaso, and Payne, as {\em generalized polyhedral complexes}, i.e., spaces glued from ``orbifold'' quotients of polyhedra with an integral affine structure factored by certain finite groups of automorphisms, \cite{ACP, CHT23}.

Let $h\colon \Gamma\to N_\R$ be a parametrized tropical curve. We define its {\em combinatorial type} $\Theta$ to be the weighted underlying graph $\mathbb G$ with ordered legs equipped with the collection of slopes $\frac{\partial h}{\partial \vec e}$ for $e \in \overline E(\G)$. An {\em isomorphism} of combinatorial types $\Theta$ and $\Theta'$ is an isomorphism of the underlying graphs respecting the order of the legs, the weight function, and the slopes of the edges and legs. Plainly, an isomorphism of parametrized tropical curves induces an isomorphism of the corresponding combinatorial types. We denote the group of automorphisms of a combinatorial type $\Theta$ by $\mathrm{Aut}(\Theta)$, and the isomorphism class of $\Theta$ by $[\Theta]$.

To a parametrized topical curve $h\colon \Gamma\to N_\R$ one can associate the point $\left(h(\mathbf{v}),\ell(\mathbf{e})\right)$ in the space $N_\mathbb R^{V(\mathbb G)}\times \mathbb R_{>0}^{E(\mathbb G)}$, and this way the set of parametrized
tropical curves of type $\Theta$ gets identified naturally with the interior $M_\Theta$ of a convex polyhedron $\overline M_\Theta\subseteq N_\mathbb R^{V(\mathbb G)}\times \mathbb R^{E(\mathbb G)}$ defined by the collection of inequalities $x(e)\ge 0$ and equalities $n({\mathfrak h}(\vec e))-n({\mathfrak t}(\vec e))=x(e)\frac{\partial h}{\partial \vec e}$ for all $\vec e\in E(\G)$; see, e.g., \cite[\S~3]{GM07}. The lattice $N^{V(\mathbb G)}\times \Z^{E(\mathbb G)}\subset N_\mathbb R^{V(\mathbb G)}\times \mathbb R^{E(\mathbb G)}$ induces a natural integral affine structure on $\overline M_\Theta$.

For each subset $E\subseteq E(\G)$, the intersection of $\overline M_\Theta$ with the locus
$$
\left\{(n(\mathbf{v}),x(\mathbf{e}))\in N_\mathbb R^{V(\mathbb G)}\times \mathbb R_{\ge 0}^{E(\mathbb G)}\; :\; x(e)=0\;\; \text{if and only if}\;\; e\in E \right\}$$
is either empty or can be identified naturally with $M_{\Theta_E}$, where $\Theta_E$ is the type of degree $\nabla$ and genus $g$, obtained from $\Theta$ by the weighted edge contraction of $E\subset E(\G)$, i.e., any connected component $\G'$ of the subgraph of $\G$ spanned by $E$ is contracted to a vertex, whose genus is defined to be $b_1(\G')+\sum_{v\in V(\G')}g(v)$. Plainly, the corresponding inclusion $\iota_{\Theta,E}\colon \overline M_{\Theta_E} \hookrightarrow \overline M_\Theta$ respects the integral affine structures.

Next, notice that an isomorphism $\alpha\colon \Theta_1\to\Theta_2$ of combinatorial types induces an isomorphism $N_\mathbb R^{V(\mathbb G_1)}\times \mathbb R^{E(\mathbb G_1)}\to N_\mathbb R^{V(\mathbb G_2)}\times \mathbb R^{E(\mathbb G_2)}$ that takes $M_{\Theta_1}$ to $M_{\Theta_2}$ and also respects the integral affine structures. Therefore, it induces an isomorphism
$\iota_\alpha\colon \overline M_{\Theta_1}\to \overline M_{\Theta_2}.$
In particular, the group $\mathrm{Aut}(\Theta)$ acts naturally on $\overline M_\Theta$.

The moduli space $M_{g, r, \nabla}^{\rm trop}$ is defined to be the colimit of the diagram, whose entries are $\overline M_\Theta$'s for all combinatorial types $\Theta$ of genus $g$ degree $\nabla$ curves having exactly $r$ contracted legs $l_1,\dotsc,l_r$; and arrows are the inclusions $\iota_{\Theta,E}$'s and the isomorphisms $\iota_\alpha$'s described above. By the construction, for each such $\Theta$, we have a finite-to-one map $M_\Theta \to M_{g, r, \nabla}^{\rm trop}$, whose image $M_\Theta/\mathrm{Aut}(\Theta)$ is denoted by $M_{[\Theta]}$.

\subsubsection{Regularity}
A parametrized tropical curve $(\Gamma,h)$ and its combinatorial type $\Theta$ are called {\em regular} if $M_\Theta$ has the expected dimension, 
\begin{equation*} \label{eq:expected_dim}
    \mathrm{expdim}(M_\Theta) := |\nabla|+r+(\rank(N)-3)\chi(\G)-{\rm ov}(\G),
\end{equation*}
where $\G$, as usual, denotes the underlying graph, $\chi(\G)$ its Euler characteristic, and ${\rm ov}(\G)$ the overvalency of $\G$, i.e., ${\rm ov}(\G):=\sum_{v\in V(\Gamma)}\max\{0,{\rm val}(v)-3\}$. 

By \cite[Proposition~2.20]{Mik05}, the expected dimension provides a lower bound on the dimension of a stratum $M_\Theta$. Note that in \emph{loc. cit.}, the expected dimension is stated to be $\mathrm{expdim}(M_\Theta) - c$, where $c$ denotes the number of edges contracted by $h$. 
This discrepancy in the formulae for the expected dimension appears because in \emph{loc. cit.} only deformations of the image $h(\Gamma)$ are taken into account. Remembering also $\Gamma$ gives one additional parameter for each contracted edge, and these parameters are independent. The following proposition is due to Mikhalkin, see \cite[Proposition~2.23]{Mik05}.

\begin{proposition}[Mikhalkin]\label{prop:mikhbound}
    Let $h\colon \Gamma\to N_\R$ be a plane tropical curve. If $h$ is an immersion away from the contracted legs and the underlying graph $\G$ is weightless and trivalent, then $(\Gamma,h)$ is regular, and $\dim(M_\Theta)=|\nabla|+r+g(\Gamma)-1$. On the other hand, if either the weight function is not identically zero or $h$ contracts some edge or ${\rm ov}(\G)>0$, then $\dim(M_\Theta)\le |\nabla|+r+g(\Gamma)-2$. 
\end{proposition}

\begin{corollary}\label{cor:gentrcur}
    In the setting and notation of the proposition, assume that $r=|\nabla|+g(\Gamma)-1$. Let $q_1,\dotsc, q_r\in N_\R$ be $r$ points in general position, i.e., $\bq=(q_1,\dotsc, q_r)$ belongs to the complement of a certain polyhedral complex of codimension-one in $N_\R^r$. If $h\colon \Gamma\to N_\R$ passes through $q_1,\dotsc, q_r$, i.e., $h(l_i)=q_i$ for all $1\le i\le r$, then the curve $h\colon \Gamma\to N_\R$ is weightless, trivalent, regular, and $h$ contracts no edges.
\end{corollary}

\subsection{The Tropicalization of Algebraic Curves}
In this section, we explain how to associate a (parametrized) tropical curve to a (parametrized) algebraic curve in a natural way. To do so, one has to work with base fields equipped with a non-archimedean valuation, e.g., the field of Puiseux series $K=\cup_{n\ge 1}\C(\!(t^{1/n})\!)$. 

Throughout the section, the ground field $K$ is the algebraic closure of a complete discretely valued field $F$ with an algebraically closed residue field. We denote the valuation on the field $K$ by $\nu\colon K\to \R\cup\{\infty\}$, the ring of integers by $K^0$, its maximal ideal by $K^{00}$, its residue field by $\tilde{K}:=K^0/K^{00}$; and similarly for $F$. 

By a {\em family of curves} we mean a {\em flat, projective} morphism of finite presentation and relative dimension one. By a collection of {\em marked points} on a family of curves we mean a collection of disjoint sections contained in the smooth locus of the family. A family of curves with marked points is {\em prestable} if its fibers have at-worst-nodal singularities; cf. \cite[\href{https://stacks.math.columbia.edu/tag/0E6T}{Tag~0E6T}]{stacks-project}. It is called {\em (semi-)stable} if so are its geometric fibers. 
Let $(C,\bp)$ be a smooth projective curve with marked points over $\Spec(K)$. By a {\em model} of $(C,\bp)$ we mean a family of curves with marked points over $\Spec(K^0)$, whose restriction to $\Spec(K)$ is $(C,\bp)$.

We fix a pair of dual lattices $M$ and $N$, and a toric variety $X$ with the lattice of characters $M$. By a {\em parametrized curve} in $X$ we mean a smooth projective curve with marked points $(C,\bp)$ and a map $f\colon  C \to X$ such that $f(C)$ does not intersect orbits of codimension greater than one, and the image of $C \setminus \bp$ under $f$ is contained in the dense torus $T \subset  X$.

\subsubsection{Tropicalization for curves}\label{subsec:tropicalization for curves}


Let $f\colon  C \to X$ be a parametrized curve, and $C^0 \to \Spec(K^0)$ a prestable model. Denote by $\widetilde{C}$ the reduction of $C$, i.e., the fiber of $C^0$ over the closed point of $\Spec(K^0)$. As usual, a point $s \in \widetilde C$ is called {\em special} if it is either a node or a marked point of $\widetilde C$. 

The {\em tropicalization} of $C$ with respect to the model $C^0$ is the tropical curve $\Gamma=(\mathbb G,\ell)$ defined as follows. The underlying graph $\G$ is the dual graph of the reduction $(\widetilde C,\tilde{\bp})$, i.e., the vertices of $\mathbb G$ correspond to irreducible components of $\widetilde{C}$, the edges -- to nodes, the legs -- to marked points, and the natural incidence relation holds: if a node belongs (resp., a marked point specializes) to an irreducible component then the corresponding edge (resp., leg) is adjacent to the corresponding vertex. In particular, a self-node of a component corresponds to a loop adjacent to the corresponding vertex. For a vertex $v$ of $\mathbb G$, its weight is defined to be the geometric genus of the corresponding component of the reduction $\widetilde{C}_v$. As for the length function, if $e\in E(\G)$ is the edge corresponding to a node $z\in \widetilde{C}$, then $\ell(e)$ is defined to be the valuation of $\lambda$, where $\lambda\in K^{00}$ is such that \'etale locally at $z$, the total space of $C^0$ is given by $xy=\lambda$. Although $\lambda$ depends on the \'etale neighborhood, its valuation does not, and hence the length function is well-defined. 

Next, we construct the parameterization $h\colon \Gamma\to N_\R$. Let $\widetilde{C}_v$ be an irreducible component of $\widetilde C$. Then, for any $m \in M$, the pullback $f^*(x^m)$ of the monomial $x^m$ is a non-zero rational function on $C^0$ since the preimage of the torus $T$ is dense in $C$. Thus, there is $\lambda_m \in K^\times$, unique up to an element invertible in $K^0$, such that $\lambda_mf^*(x^m)$ is an invertible function at the generic point of $\widetilde C_v$. The function $h(v)$, associating to $m \in M$ the valuation $\nu(\lambda_m)$, is clearly linear, and hence $h(v)\in N_\R$. The parameterization $h\colon \Gamma\to N_\R$ is defined to be the unique piecewise affine function with values $h(v)$ at the vertices of $\Gamma$, whose slopes along the legs satisfy the following: for any leg $l$ and $m\in M$ we have $\frac{\partial h}{\partial \vec l}(m)=-\mathrm{ord}_{p} f^*(x^m)$, where $p$ is the marked point corresponding to $l$. In particular, $h$ contracts $l$ if and only if $f(p)\in T$.

\begin{lemma}
    The map $h\colon \Gamma\to N_\R$ is a parametrized tropical curve, i.e., the slopes of $h$ are integral and satisfy the balancing condition at any vertex. 
\end{lemma}

The proof (stated in a slightly different language) can be found in \cite[Lemma~2.23]{Tyo12}. The curve $\Gamma$ (resp. $h\colon \Gamma\to N_\R$) is called the \emph{tropicalization} of $C$ (resp. $f\colon  C \to X$) with respect to the model $C^0$. Plainly, the tropical curve $\Gamma$ is independent of the parameterization, and depends only on $C^0$. If the curve $(C,\bp)$ is stable and $C^0$ is the stable model, then the corresponding tropicalization is called simply {\em the tropicalization} of $(C,\bp)$ (resp. $f\colon  C \to X$), and is denoted by ${\rm trop}(C)$ (resp. $h\colon {\rm trop}(C)\to N_\R$). Plainly, ${\rm trop}(C)$ is stable.

Denote by ${\rm Val}\colon T(K)\to N_\R$ the coordinatewise valuation map from the dense torus $T\subset X$ to $N_\R$ given by ${\rm Val}(p)(m):=\nu(x^m(p))$. And set ${\rm Trop}:=-{\rm Val}$. The map ${\rm Trop}$ is often called the {\em tropicalization map} of the algebraic torus $T$, and for any subvariety $Y\subset T$, the closure of ${\rm Trop}(Y)$ in $N_\R$ is called the {\em non-Archimedean amoeba} of $Y$.

\begin{exercise}
    Let $(C,\bp;f\colon  C \to X)$ be a parametrized curve, $h\colon \Gamma\to N_\R$ its tropicalization, and $\nabla$ the degree of $(\Gamma, h)$. Let $l$ be the leg corresponding to a marked point $p$. Show that 
    \begin{itemize}
        \item If $h$ contracts $l$, then $h(l)(m)=-\nu(x^m(f(p)))$ for all $m\in M$; 
        \item The image of $h$ is the non-Archimedean amoeba of $f(C)\cap T$; 
        \item If $h$ does not contract $l$, then $f$ maps $p$ to a codimension-one orbit. Furthermore, if $O$ is the orbit containing $f(p)$ and $n$ is the primitive generator of the ray in the fan of $X$ corresponding to $O$, then the slope of $h$ along $l$ is given by $\frac{\partial h}{\partial \vec l}=-\ord_p(f^*(O))n$. In particular, if $f$ is transversal to the boundary divisor, then the degree $\nabla$ is reduced;
        \item If $X=X_\Delta$ is a projective toric surface and $C$ the normalization of a curve in the linear system $|\ko(\Delta)|$, then the degree $\nabla$ is dual to $\Delta$.
    \end{itemize}
\end{exercise}

\subsubsection{A baby example}
Let $C\subset\P^2$ be the line given by the homogeneous equation $x+\mu y=z$, for some $\mu\in K^{00}$, and let $f$ be the immersion of the line $C$ in $\P^2$. Let us mark the following points on the curve $C$: $p_1:=[1-\mu:1:1], p_2:=[1:0:1], p_3:=[-\mu:1:0], p_4:=[0:1:\mu]$. The homogeneous equation $x+\mu y=z$ defines a trivial model of the curve $C$ over the ring of integers $K^0$, whose reduction has a unique component $L$ (the line in the complex plane given by the equation $x=z$). Plainly, $p_3$ and $p_4$ specialize to the same point $[0:1:0]$. Hence this integral model is not even a prestable model of the marked curve $(C,\bp)$. 

In the affine coordinates, the trivial model is given by $$\Spec\,K\left[\frac{x}{y},\frac{z}{y}\right]/\left(\frac{x}{y}-\frac{z}{y}+\mu\right)\simeq \Spec\,K\left[\frac{z}{y}\right]$$ and the marked points $p_3, p_4$ are given by $\frac{z}{y}=0$ and $\frac{z}{y}=\mu$, respectively. After blowing up this model with respect to the ideal generated by $\frac{z}{y}$ and $\mu$, the points $p_3$ and $p_4$ get separated and specialize to two distinct points of the exceptional divisor $E$, while the points $p_1$ and $p_2$ specialize to two distinct points of the strict transform of $\tilde{L}$. Furthermore, all four specializations are different from the node of the reduction. Thus, we have constructed the stable model of $(C,\bp)$ and can now describe its tropicalization. The graph $\G$ consists of two vertices, joined by a single edge, and four ends, with the first two ends being attached to one of the vertices and the last two to another. Finally, the length function $\ell$ obtains value $\nu(\mu)$ on the unique edge since (\'etale) locally the model is given by $x_1x_2=\mu$ (the local coordinates here are $x_1:=\frac{z}{y}$ and $x_2:=\frac{\mu y}{z}$). See Figure~\ref{fig:trop} for an illustration.

It remains to describe the parametrization map $h\colon {\rm trop}(C)\to \R^2$. Let $v$ be the vertex corresponding to the component $\tilde{L}$ and $w$ to $E$. Since the monomial functions $x^{(1,0)}=\frac{x}{z}$ and $x^{(0,1)}=\frac{y}{z}$ are invertible at the generic point of $\tilde{L}$, we have $h(v)=(0,0)$. Similarly, since $\frac{\mu y}{z}$ is a coordinate function on the exceptional divisor $E$ it is invertible at the generic point of $E$ and hence so is $\frac{x}{z}=1-\frac{\mu y}{z}$. Thus, $h(w)=(0,\nu(\mu))$. It remains to calculate the slopes of $h$ along the legs. It follows from the definition that $\frac{\partial h}{\partial \vec l_1}=(0,0)$, i.e., $h$ contracts $l_1$, $\frac{\partial h}{\partial \vec l_2}=-(0,1)$, $\frac{\partial h}{\partial \vec l_3}=-(-1,-1)=(1,1)$, and $\frac{\partial h}{\partial \vec l_4}=-(1,0)$. Now it is evident that all slopes of $h$ are integral, and that the balancing condition is satisfied at any vertex.

\tikzset{every picture/.style={line width=0.75pt}} 
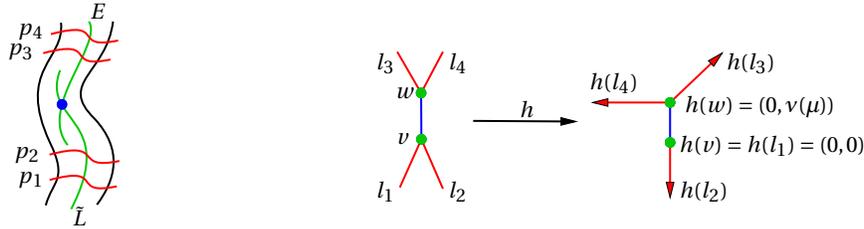
\begin{figure}[ht]
\centering
\begin{tikzpicture}[x=0.48pt,y=0.48pt,yscale=-1,xscale=1]
\draw [line width=0.75]    (139.2,21.6) .. controls (134.5,62.1) and (95.2,69.6) .. (124.5,98.1) .. controls (147.2,122.6) and (141.5,146.1) .. (126.5,166.1) ;
\draw [color={rgb, 255:red, 5; green, 200; blue, 5}  ,draw opacity=1 ][line width=0.75]    (120.2,21.6) .. controls (131.5,40.1) and (79.5,89.1) .. (103.5,119.1) ;
\draw [color={rgb, 255:red, 5; green, 200; blue, 5}  ,draw opacity=1 ][line width=0.75]    (97.5,60.1) .. controls (81.8,101.6) and (143.5,97.1) .. (106.2,169.6) ;
\draw [line width=0.75]    (98.2,22.6) .. controls (101.5,52.1) and (62.5,59.1) .. (89.5,117.1) .. controls (103.8,146.6) and (91.5,147.1) .. (86.2,166.6) ;
\draw [color={rgb, 255:red, 255; green, 0; blue, 0 }  ,draw opacity=1 ][line width=0.75]    (90.2,31.6) .. controls (135.2,19.6) and (105.2,44.6) .. (143.2,36.6) ;
\draw [color={rgb, 255:red, 255; green, 0; blue, 0 }  ,draw opacity=1 ][line width=0.75]    (84.2,46.6) .. controls (129.2,34.6) and (99.2,59.6) .. (137.2,51.6) ;
\draw [color={rgb, 255:red, 255; green, 0; blue, 0 }  ,draw opacity=1 ][line width=0.75]    (89.2,126.6) .. controls (134.2,114.6) and (106.2,139.6) .. (144.2,131.6) ;
\draw [color={rgb, 255:red, 255; green, 0; blue, 0 }  ,draw opacity=1 ][line width=0.75]    (89.2,146.6) .. controls (134.2,134.6) and (104.2,159.6) .. (142.2,151.6) ;
\draw [color={rgb, 255:red, 0; green, 0; blue, 255 }  ,draw opacity=1 ][line width=0.75]    (99,87) ;
\draw [shift={(99,87)}, rotate = 0] [color={rgb, 255:red, 0; green, 0; blue, 255 }  ,draw opacity=1 ][fill={rgb, 255:red, 0; green, 0; blue, 255 }  ,fill opacity=1 ][line width=0.75]      (0, 0) circle [x radius= 3.35, y radius= 3.35]   ;
\draw [shift={(99,87)}, rotate = 0] [color={rgb, 255:red, 0; green, 0; blue, 255 }  ,draw opacity=1 ][fill={rgb, 255:red, 0; green, 0; blue, 255 }  ,fill opacity=1 ][line width=0.75]      (0, 0) circle [x radius= 3.35, y radius= 3.35]   ;
\draw [color={rgb, 255:red, 255; green, 0; blue, 0 }  ,draw opacity=1 ][line width=0.75]    (363.2,45.5) -- (382,78) ;
\draw [color={rgb, 255:red, 255; green, 0; blue, 0 }  ,draw opacity=1 ][line width=0.75]    (399.2,45.5) -- (382,78) ;
\draw [color={rgb, 255:red, 255; green, 0; blue, 0 }  ,draw opacity=1 ][line width=0.75]    (382.2,114.5) -- (365.2,152.5) ;
\draw [color={rgb, 255:red, 255; green, 0; blue, 0 }  ,draw opacity=1 ][line width=0.75]    (382.2,114.5) -- (399.2,153.5) ;
\draw [line width=0.75]    (423,100) -- (502.2,100.49) ;
\draw [shift={(504.2,100.5)}, rotate = 180.35] [fill={rgb, 255:red, 0; green, 0; blue, 0 }  ][line width=0.08]  [draw opacity=0] (12,-3) -- (0,0) -- (12,3) -- cycle    ;
\draw [color={rgb, 255:red, 255; green, 0; blue, 0 }  ,draw opacity=1 ][line width=0.75]    (578.2,85.5) -- (617.75,47.88) ;
\draw [shift={(619.2,46.5)}, rotate = 136.43] [fill={rgb, 255:red, 255; green, 0; blue, 0 }  ,fill opacity=1 ][line width=0.08]  [draw opacity=0] (12,-3) -- (0,0) -- (12,3) -- cycle    ;
\draw [color={rgb, 255:red, 255; green, 0; blue, 0 }  ,draw opacity=1 ][line width=0.75]    (519.2,85.5) -- (578.2,85.5) ;
\draw [shift={(517.2,85.5)}, rotate = 0] [fill={rgb, 255:red, 255; green, 0; blue, 0 }  ,fill opacity=1 ][line width=0.08]  [draw opacity=0] (12,-3) -- (0,0) -- (12,3) -- cycle    ;
\draw [color={rgb, 255:red, 255; green, 0; blue, 0 }  ,draw opacity=1 ][line width=0.75]    (578.2,158.5) -- (578.2,117) ;
\draw [shift={(578.2,160.5)}, rotate = 270] [fill={rgb, 255:red, 255; green, 0; blue, 0 }  ,fill opacity=1 ][line width=0.08]  [draw opacity=0] (12,-3) -- (0,0) -- (12,3) -- cycle    ;
\draw [color={rgb, 255:red, 0; green, 0; blue, 255 }  ,draw opacity=1 ][line width=0.75]    (578.2,117) -- (578.2,85.5) ;
\draw [shift={(578.2,85.5)}, rotate = 270] [color={rgb, 255:red, 0; green, 0; blue, 255 }  ,draw opacity=1 ][fill={rgb, 255:red, 0; green, 0; blue, 255 }  ,fill opacity=1 ][line width=0.75]      (0, 0) circle [x radius= 3.35, y radius= 3.35]   ;
\draw [shift={(578.2,117)}, rotate = 270] [color={rgb, 255:red, 0; green, 0; blue, 255 }  ,draw opacity=1 ][fill={rgb, 255:red, 0; green, 0; blue, 255 }  ,fill opacity=1 ][line width=0.75]      (0, 0) circle [x radius= 3.35, y radius= 3.35]   ;
\draw [color={rgb, 255:red, 0; green, 0; blue, 255 }  ,draw opacity=1 ][line width=0.75]    (382,78) -- (382.2,114.5) ;
\draw [shift={(382.2,114.5)}, rotate = 89.69] [color={rgb, 255:red, 0; green, 0; blue, 255 }  ,draw opacity=1 ][fill={rgb, 255:red, 0; green, 0; blue, 255 }  ,fill opacity=1 ][line width=0.75]      (0, 0) circle [x radius= 3.35, y radius= 3.35]   ;
\draw [shift={(382,78)}, rotate = 89.69] [color={rgb, 255:red, 0; green, 0; blue, 255 }  ,draw opacity=1 ][fill={rgb, 255:red, 0; green, 0; blue, 255 }  ,fill opacity=1 ][line width=0.75]      (0, 0) circle [x radius= 3.35, y radius= 3.35]   ;
\draw [line width=0.75]    (382,78) ;
\draw [shift={(382,78)}, rotate = 0] [color={rgb, 255:red, 5; green, 200; blue, 5 }  ][fill={rgb, 255:red, 5; green, 200; blue, 5 }  ][line width=0.75]      (0, 0) circle [x radius= 3.35, y radius= 3.35]   ;
\draw [shift={(382,78)}, rotate = 0] [color={rgb, 255:red, 5; green, 200; blue, 5 }  ][fill={rgb, 255:red, 5; green, 200; blue, 5 }  ][line width=0.75]      (0, 0) circle [x radius= 3.35, y radius= 3.35]   ;
\draw [line width=0.75]    (382.2,114.5) ;
\draw [shift={(382.2,114.5)}, rotate = 0] [color={rgb, 255:red, 5; green, 200; blue, 5 }  ][fill={rgb, 255:red, 5; green, 200; blue, 5 }  ][line width=0.75]      (0, 0) circle [x radius= 3.35, y radius= 3.35]   ;
\draw [shift={(382.2,114.5)}, rotate = 0] [color={rgb, 255:red, 5; green, 200; blue, 5 }  ][fill={rgb, 255:red, 5; green, 200; blue, 5 }  ][line width=0.75]      (0, 0) circle [x radius= 3.35, y radius= 3.35]   ;
\draw [line width=0.75]    (578.2,117) ;
\draw [shift={(578.2,117)}, rotate = 0] [color={rgb, 255:red, 5; green, 200; blue, 5 }  ][fill={rgb, 255:red, 5; green, 200; blue, 5 }  ][line width=0.75]      (0, 0) circle [x radius= 3.35, y radius= 3.35]   ;
\draw [shift={(578.2,117)}, rotate = 0] [color={rgb, 255:red, 5; green, 200; blue, 5 }  ][fill={rgb, 255:red, 5; green, 200; blue, 5 }  ][line width=0.75]      (0, 0) circle [x radius= 3.35, y radius= 3.35]   ;
\draw [line width=0.75]    (578.2,85.5) ;
\draw [shift={(578.2,85.5)}, rotate = 0] [color={rgb, 255:red, 5; green, 200; blue, 5 }  ][fill={rgb, 255:red, 5; green, 200; blue, 5 }  ][line width=0.75]      (0, 0) circle [x radius= 3.35, y radius= 3.35]   ;
\draw [shift={(578.2,85.5)}, rotate = 0] [color={rgb, 255:red, 5; green, 200; blue, 5 }  ][fill={rgb, 255:red, 5; green, 200; blue, 5 }  ][line width=0.75]      (0, 0) circle [x radius= 3.35, y radius= 3.35]   ;

\draw (344,44) node [anchor=north west][inner sep=0.75pt]  [font=\footnotesize]  {$l_{3}$};
\draw (401,44) node [anchor=north west][inner sep=0.75pt]  [font=\footnotesize]  {$l_{4}$};
\draw (344,145) node [anchor=north west][inner sep=0.75pt]  [font=\footnotesize]  {$l_{1}$};
\draw (401,145) node [anchor=north west][inner sep=0.75pt]  [font=\footnotesize]  {$l_{2}$};
\draw (457,82) node [anchor=north west][inner sep=0.75pt]  [font=\footnotesize]  {$h$};
\draw (359,71) node [anchor=north west][inner sep=0.75pt]  [font=\footnotesize]  {$w$};
\draw (360,107) node [anchor=north west][inner sep=0.75pt]  [font=\footnotesize]  {$v$};
\draw (583,109) node [anchor=north west][inner sep=0.75pt]  [font=\footnotesize,rotate=-0.15]  {$h(v)=h(l_1)=(0,0)$};
\draw (587,78) node [anchor=north west][inner sep=0.75pt]  [font=\footnotesize]  {$h(w)=( 0,\nu ( \mu ))$};
\draw (583,145) node [anchor=north west][inner sep=0.75pt]  [font=\footnotesize]  {$h( l_{2})$};
\draw (512,60) node [anchor=north west][inner sep=0.75pt]  [font=\footnotesize]  {$h( l_{4})$};
\draw (620,44) node [anchor=north west][inner sep=0.75pt]  [font=\footnotesize]  {$h( l_{3})$};
\draw (54,37) node [anchor=north west][inner sep=0.75pt]  [font=\footnotesize]  {$p_{3}$};
\draw (62,138) node [anchor=north west][inner sep=0.75pt]  [font=\footnotesize]  {$p_{1}$};
\draw (59,118) node [anchor=north west][inner sep=0.75pt]  [font=\footnotesize]  {$p_{2}$};
\draw (62,21) node [anchor=north west][inner sep=0.75pt]  [font=\footnotesize]  {$p_{4}$};
\draw (105,165) node [anchor=north west][inner sep=0.75pt]  [font=\footnotesize]  {$\tilde{L}$};
\draw (119,5) node [anchor=north west][inner sep=0.75pt]  [font=\footnotesize]  {$E$};

\end{tikzpicture}	

\caption{On the left, we have a cartoon of the stable model of $(C,\bp)$, and on the right, its tropicalization.}\label{fig:trop}
\end{figure}

\begin{exercise}
    Let $\nabla$ be the degree of the parametrized tropical curve $h\colon \Gamma\to \R^2$ constructed above. Describe the moduli space $M_{0, 4, \nabla}^{\rm trop}$ explicitly, and show that the curve $h\colon \Gamma\to \R^2$ is regular. 
\end{exercise}

\subsection{The Multiplicity of a Tropical Curve}
In this subsection we define Mikhalkin's multiplicities of plane tropical curves and state Mikhalkin's Correspondence Theorem, \cite[Theorem~1]{Mik05}.

Let $h\colon \Gamma\to N_\R$ be a parametrized tropical curve, $v$ a vertex of valency $3$, and $\vec{e}_1,\vec{e}_2\in{\rm Star}(v)$ two distinct oriented edges. Denote by $\mathfrak{A}(v)$ the area of the parallelogram spanned by the slopes $\frac{\partial h}{\partial \vec{e}_1}$ and $\frac{\partial h}{\partial \vec{e}_2}$, and notice that this area is independent of the choice of $\vec{e}_1$ and $\vec{e}_2$ by the balancing condition.

\begin{definition}
    The {\em Mikhalkin's multiplicity} of a weightless and trivalent plane tropical curve $h\colon \Gamma\to N_\R$ is the product $\prod \mathfrak{A}(v)$ running over all vertices of $\Gamma$ not adjacent to contracted legs. We denote this multiplicity by $\mathfrak{mult}(\Gamma,h)$.
\end{definition}

\begin{theorem}[Mikhalkin's Correspondence Theorem]
Assume that the residue filed of $K$ has characteristic zero. Let $\Delta\subset M_\R$ be a lattice polygon, $g$ an integer, and $\nabla$ a reduced tropical degree dual to $\Delta$. Set $r:=|\nabla|+g-1=|\partial\Delta\cap M|+g-1$. Let $t_1,\dotsc, t_r\in T\subset X$ be points in general position, such that their tropicalizations $q_i:={\rm Trop}(t_i)$ are in tropically general position in $N_\R$. Let $h\colon \Gamma\to N_\R$ be a parametrized tropical curve of degree $\nabla$, genus $g$, and with exactly $r$ contracted legs $l_1,\dotsc, l_r$ passing through $q_1,\dotsc, q_r$. Then the number of parametrized algebraic curves $(C,\bp; f\colon C\to X)$ passing through the $t_i$'s and tropicalizing to $(\Gamma, h)$ is equal to Mikhalkin's multiplicity $\mathfrak{mult}(\Gamma,h)$.
\end{theorem}

\begin{remark}
    (1) By Corollary~\ref{cor:gentrcur}, the tropical curves appearing in the Correspondence Theorem are weightless and trivalent. Therefore, their multiplicities are well-defined. 

    (2) An algebraic proof of the theorem (stated in a slightly different language) can be found in \cite[Theorem~6.2]{Tyo12}. In fact, the proof in {\em loc. cit.} is valid also if the residue characteristic is big enough, i.e., does not divide the multiplicity $\mathfrak{mult}(\Gamma,h)$. See also \cite{Shu05}.
    
    (3) Mikhalkin's enumerative formula (Theorem~\ref{thm:MikFor}) follows from Mikhalkin's Correspondence Theorem, but one shall be slightly careful. Namely, in Theorem~\ref{thm:MikFor}, one counts the isomorphism classes of parametrized tropical curves with {\em unordered} non-contracted legs since the intersection of curves parametrized by the Severi variety with the boundary divisor comes unordered. Therefore, after counting the isomorphism classes of parametrized tropical curves with Mikhalkin's multiplicities, one has to divide the result by the number of permutations $\sigma\in {\mathfrak S}_{|\nabla|}$ satisfying $\frac{\partial h}{\partial \vec{l}_{r+j}}=\frac{\partial h}{\partial \vec{l}_{r+\sigma(j)}}$ for all $j$. Since $\nabla$ is reduced, the latter number is nothing but $\prod_i (k_i!)$, where $k_i=|\partial_i\Delta\cap M|$. 
\end{remark}

\subsection{Tropical Methods  -- Beyond the Enumeration}

Following the pioneering work of Mikhalkin \cite{Mik05}, tropical techniques have been used to get a characteristic independent proof of the upper bound on the dimension of Severi varieties on toric surfaces and to resolve the Severi problem in positive characteristic, \cite{CHT23}, to prove the irreducibility of Hurwitz spaces in small characteristics, \cite{CHT24}, and to study the adjacencies between the Severi varieties of different genera in a given linear system on toric surfaces, \cite{CHT22t}.

\subsubsection{Dimension bounds}\label{sect:proofdimest}
In this section, we show how to use the tropical techniques introduced in the previous section to prove that $\dim(V_{g,\Delta})\le|\partial\Delta\cap M|+g-1$. Assume to the contrary that $\dim(V_{g,\Delta})\ge|\partial\Delta\cap M|+g$. Set $r:=|\partial\Delta\cap M|+g$, and let $t_1,\dotsc, t_r\in T\subset X_\Delta$ be points in general position, such that their tropicalizations $q_i:={\rm Trop}(t_i)$ are in tropically general position in $N_\R$. For any $i$, let $H_i\subset |\ko(\Delta)|$ be the hyperplane of curves passing through the point $t_i$. Then, by dimension reasons, the closure of $V_{g,\Delta}$ intersects $\cap_i H_i$ non-trivially in the projective space $|\ko(\Delta)|$. And since the $t_i$'s are in general position, $\left(\cap_i H_i\right)\cap V_{g,\Delta}\ne\emptyset$. Thus, there exists a parametrized curve $(C,\bp;f\colon C\to X_\Delta)$ of genus $g$ that passes through the points $t_1,\dotsc, t_r$, i.e., $f(p_i)=t_i$ for all $1\le i\le r$.

Consider the tropicalization $h:\Gamma\to N_\R$ of such a curve $(C,\bp;f\colon C\to X_\Delta)$. By the construction, $(\Gamma, h)\in M_{g, r, \nabla}^{\rm trop}$, where $\nabla$ is a (reduced) tropical degree dual to $\Delta$. Furthermore, $h(l_i)=q_i$ for any $1\le i\le r$. Now, since the points $q_1,\dotsc, q_r\in N_\R$ are in tropically general position, it follows that 
$$\dim M_{g, r, \nabla}\ge r\cdot \rank(N)=2r=|\partial\Delta\cap M|+r+g>|\nabla|+r+g(\Gamma)-1,$$
which is a contradiction by Proposition~\ref{prop:mikhbound}.

\subsubsection{Applications to the irreducibility problems}

There are currently two tropical approaches to the proof of the irreducibility of Severi varieties developed in \cite{CHT23, CHT24}. Both approaches are technically involved. Thus, let us briefly discuss the ideas behind these approaches without going into technicalities. 

In \cite{CHT23}, a tropical approach to the Severi problem is developed and implemented in the classical case of the degree-$d$ curves on the projective plane. The general framework follows the original strategy of Severi. Namely, one first shows that any irreducible component of the Severi variety $V_{g,d}$ contains the variety $V_{1-d,d}$ in its closure, and then uses monodromy type arguments to show that there exists a unique component of $V_{g,d}^{\rm irr}$ containing $V_{1-d,d}$ in its closure. 

The second part of the argument is relatively easy and works in arbitrary characteristic. The first part of the argument, however, is much more involved. In Harris' solution of the Severi problem, this part uses deformation-theoretic arguments, computation of first-order deformations, and the study of the versal deformation space of tacnodes -- all these are characteristic sensitive, and we do not know how to generalize Harris' arguments to the case of positive characteristic. In \cite{CHT23}, this part of the argument is replaced with tropical techniques. 

The main technical tool developed in \cite{CHT23} is the notion of the {\em tropicalization of one-parameter families of parametrized tropical curves} and the proof of geometric properties of the induced map $\alpha$ from the tropicalization of the base to the moduli space $M_{g, r, \nabla}^{\rm trop}$. It turns out that the map $\alpha$ is piecewise integral affine, and in good cases, it satisfies certain {\em balancing properties}. Roughly speaking, what these properties allow one to do is to control the intersection of the image of $\alpha$ with the maximal-dimensional strata in $M_{g, r, \nabla}^{\rm trop}$. In particular, one can show that for a maximal-dimensional stratum, the image of $\alpha$ is either disjoint from it or intersects it along a straight interval, whose boundary does not belong to the stratum. Furthermore, if this interval intersects the boundary along a codimension-one stratum $M_{[\Theta]}$, then there are two more maximal-dimensional strata containing $M_{[\Theta]}$ in their closures, and the image of $\alpha$ intersects both of them. 

To prove that the closure $\overline V\subseteq |\ko_{\P^2}(d)|$ of any irreducible component $V\subseteq V_{g,d}$ contains $V_{1-d,d}$, one proceeds by induction and shows that $\overline V$ contains an irreducible component of $V_{g-1,d}$, and therefore ultimately a component of $V_{1-d,d}$. Since the variety $V_{1-d,d}$ parametrizes unions of $d$ distinct lines, it is evidently irreducible, and therefore, $V_{1-d,d}\subset \overline V$. 

To prove that $\overline V$ contains an irreducible component of $V_{g-1,d}$, one considers the intersection $Z$ of $V$ with the space of curves passing through $3d+g-2$ points in general position, chosen such that their tropicalizations are in tropically general position. One then considers the tropicalization of the tautological family over the curve $Z$ and uses the general properties of tropicalizations discussed above to reduce the degeneration assertion to a combinatorial game with tropical curves and their moduli. 

The goal of the combinatorial game is to prove that the image of $\alpha$ must intersect a maximal-dimensional stratum parameterizing tropical curves having a contracted edge of varying length. Indeed, if this is the case, then ${\rm trop}(Z^\nu)$ contains a leg parameterizing tropical curves for which the length of the contracted edge is growing to infinity. If $[C]\in Z$ is the marked point corresponding to such a leg, then it is not difficult to check that $C$ is necessarily reduced and has geometric genus $g-1$. Finally, since $C$ passes through $3d+g-2$  points in general position, it follows by dimension reasoning that $\overline V$ contains a component of the Severi variety $V_{g-1,d}$ as needed.

\begin{remark}
Let $\Delta\subset\R^2$ be an $h$-transverse polygon, i.e., a lattice polygon, whose intersections with the horizontal lines of integral height are either empty or intervals with integral boundary points. Let $0\le g\le |\Delta^\circ\cap\Z^2|$ be an integer, and consider the Severi variety $V_{g,\Delta}^{\rm irr}$. Then one can use the balancing properties discussed above together with some subtler balancing properties of the map $\alpha$ in order to show that any irreducible component of $V_{g,\Delta}^{\rm irr}$ contains $V_{0,\Delta}^{\rm irr}$ in its closure provided that the characteristic of the ground field is either zero or big enough, see \cite{CHT22t}. It follows that the closure of $V_{g,\Delta}^{\rm irr}$ contains at least one irreducible component of $V_{g',\Delta}^{\rm irr}$ for any $0\le g'\le g$. Furthermore, for each prime $p$, in {\em loc. cit.} there is an example of an $h$-transverse polygon for which there exists a component of $V_{1,\Delta}^{\rm irr}$ whose closure is disjoint from $V_{0,\Delta}^{\rm irr}$.
\end{remark}

A different tropical approach to the irreducibility problem is developed in \cite{CHT24} and implemented in the case of Severi varieties on polarized toric surfaces associated to almost all $h$-transverse polygons, e.g., Hirzebruch surfaces, toric del-Pezzo surfaces, etc. As an application, the irreducibility of Hurwitz spaces\footnote{Recall that Hurwitz space $H_{g,d}$ is the moduli space of genus-$g$ degree-$d$ simple branched tame coverings of the projective line $\P^1$, where by simple we mean that any point has at least $d-1$ preimages.} $H_{g,d}$ is deduced in any characteristic, completing Fulton's result from 1969 that proves the irreducibility of $H_{g,d}$ under the assumption that the characteristic of the ground field is at least $d+1$, see \cite{Ful69}.  

Set $r:=\dim(V_{g,\Delta}^{\rm irr})$. Pick a trivalent parametrized tropical curve $(\Gamma,h)$ of genus $g$, reduced degree $\nabla$ dual to $\Delta$, and with exactly $r$ contracted ends for which Mikhalkin's multiplicity is one and such that $q_i:=h(l_i)$ for $1\le i\le r$ are in tropically general position. Pick some points $t_1,\dotsc, t_r\in T\subset X_\Delta$ such that ${\rm Trop}(t_i)=q_i$ for all $1\le i\le r$. We can now describe the general strategy of \cite{CHT24} to the proof of the irreducibility of $V_{g,\Delta}^{\rm irr}$ without going into technical details. 

The main idea is to use the properties of the tropicalizations of families of curves and the combinatorial properties of the moduli space $M_{g, r, \nabla}^{\rm trop}$ to prove that there exists a parametrized curve $(C, \bp; f\colon  C\to X_\Delta)$ that passes through the points $t_1,\dotsc, t_r$, tropicalizes to $(\Gamma,h)$, and such that $[f(C)]$ belongs to a given irreducible component $V\subseteq V_{g,\Delta}^{\rm irr}$. The multiplicity-one property of $(\Gamma,h)$ is then used to show that there exists a unique irreducible component containing $[f(C)]$. Notice that unlike Severi's original approach implemented by Harris and followed in \cite{CHT23}, this new approach contains neither degeneration arguments nor monodromy-type arguments.

\egroup
				  
\bibliographystyle{amsalpha}
\providecommand{\bysame}{\leavevmode\hbox to3em{\hrulefill}\thinspace}
\providecommand{\MR}{\relax\ifhmode\unskip\space\fi MR }
\providecommand{\MRhref}[2]{%
  \href{http://www.ams.org/mathscinet-getitem?mr=#1}{#2}
}
\providecommand{\href}[2]{#2}

\end{document}